\documentclass[12pt,reqno]{amsart}
\usepackage{amsmath}
\usepackage{amssymb}

\thanks{
J.S.'s research supported in part by
the Swedish Natural Science Research
Council and the G\"{o}ran Gustafsson Foundation for Research in 
Natural Sciences and Medicine.}

\author[Jeffrey E. Steif]{Jeffrey E. Steif}
\address[Jeffrey E. Steif]{Department of Mathematics, Chalmers University
of Technology, S-412 96 Gothenburg, Sweden}
\email[Jeffrey E. Steif]{steif@math.chalmers.se}

\author[Aidan Sudbury]{Aidan Sudbury}
\address[Aidan Sudbury]
{ School of Mathematical Sciences, Monash University, Vic 3800,
Australia}
\email[Aidan Sudbury]{aidan.sudbury@sci.monash.edu.au}

\usepackage{latexsym}
\usepackage{amsfonts}
\usepackage{amsmath}
\usepackage{amsthm}
\usepackage{epsfig}

\usepackage{graphicx}
\usepackage{amsmath}
\newtheorem{theorem}{Theorem}

\newtheorem{conjecture}[theorem]{Conjecture}

\newtheorem{lemma}[theorem]{Lemma}

\begin{document}

\title[Poisoning in a catalytic model]
{Some results for poisoning in a catalytic model}

\begin{abstract}
We obtain new results concerning poisoning/nonpoisoning in
a catalytic model which has previously been introduced and studied. 
We show that poisoning can occur even when the arrival rate of one gas 
is smaller than the sum of the arrival rates of the other gases, 
and that poisoning does not occur when all gases have equal
arrival rates.

\end{abstract}

\date{April 11, 2006}

\maketitle
\noindent



\section{\textbf{Introduction}}
Grannan and Swindle (1991) introduced a collection of interacting
particle systems to model catalytic surfaces. In these models $n$
types of gases represented by states $\{1,2,...,n\}$ fall on sites of
the integer lattice $Z^d$. A vacant site is represented by state
$0$. Molecules of gas $i$ fall upon vacant sites at rate $p_i$,
where $\sum p_i=1$. No two different gases can occupy adjacent
sites and so if a molecule of type $i$ falls upon a vacant site
adjacent to a molecule of type $j,j\neq i$, then the two gases
react and both sites are left vacant. If there are several
adjacent sites with molecules different from the arriving molecule, one of them
is selected uniformly at random to react with the arriving molecule.

Grannan and Swindle (1991) defined "poisoning" as the configuration converging
a.s.\ as time goes to infinity. In such a case, the limit is necessary the 
``all $i$ configuration'' for some $i\in \{1,2,...,n\}$.
Let $\delta_i$ be the point mass at the configuration where everything 
is in state $i$ and define {\sl coexistence} as the existence
of a stationary distribution which is not a 
convex combination of the $\{\delta_i\}$'s.
Heuristically, poisoning implies noncoexistence 
(but this statement is not precise). 
Grannan and Swindle (1991) proved the following result. (They 
proved somewhat more than what we state.)

\begin{theorem} 
(Grannan and Swindle)
(i). If $n=2$, $d$ is arbitrary and $p_1\neq 1/2$, then coexistence does 
not occur. \\
(ii). There exists $\epsilon_0>0$ such that for any $n$ and $d$, if
$p_1>1-\epsilon_0$, then for any initial state having infinitely many 0's,
a.s. poisoning occurs with the limiting configuration being all 1's;
i.e., $\eta_t \rightarrow 1$, the configuration consisting of all 1's. \\
(iii). For any $n$ and $d$, if $p_1 < \frac{1}{2d\lambda_c(d)+1}$ (where
$\lambda_c(d)$ is the critical value for the contact process in $d$--dimensions),
then, for any initial state having infinitely many 0's, a.s.\ the system does not 
converge to the all 1 configuration. \\
(iv). For every $d$, there exists $N_d$ such that if $n \ge N_d$, there exist
$p_1,\ldots,p_n$ so that coexistence occurs. For $d=1$, one can have coexistence
with $n=5$.
\end{theorem}

(i) was proved using a type of ``energy'' argument, (ii) was quite involved and
(iii) and (iv) were proved by arguments involving stochastic comparison 
with a supercritical contact process. 
It is not known whether coexistence can occur for less than 5 gases in 1 dimension.

Mountford and Sudbury (1992) strengthened part (ii) (using a submartingale argument)
of the above result by showing:
\begin{theorem} \label{MS}
(Mountford and Sudbury)
 Let $n$ and $d$ be arbitrary and $p_1> 1/2$. Then,
if the initial state,$\eta_0$, has infinitely many
1's or 0's, then poisoning occurs with $\eta_t \rightarrow 1$ a.s.
\end{theorem}

A question left open is whether it is possible for
a gas to poison a surface with its arrival rate being less than 1/2.
Our first result proves this to be the case.

\begin{theorem} \label{first theorem}
If $d=1$ and $n=4$, if the first gas has rate $0.47$
and the other three have rates $0.53/3$, then the first gas
will poison the surface with probability 1 if the initial state of
the surface has an infinite number of 0's or 1's.
\end{theorem}

It is trivial to observe that if poisoning occurs from every initial 
state, then there cannot be coexistence. 
However, one should not think of these things as being synonymous as we will see.
In addition, the proof of  Theorem \ref{MS}
shows that for any $d$ and $n$, if $p_1>1/2$, then coexistence does not occur.
For the interesting case $p_1=p_2$, Mountford (1992) proved the following result.
(What we state is just a special case of what he proved.)

\begin{theorem}
(Mountford)
For $d=1$ and $n=2$, the only translation invariant stationary distributions are
$\{p\delta_1+(1-p)\delta_2\}_{p\in [0,1]}$.
\end{theorem}

This {\sl strongly} suggests that coexistence does not occur here. Our second
result says that this noncoexistence result is not due to poisoning.

\begin{theorem} \label{second theorem}
Let $n$ and $d$ be arbitrary. If there is $i\neq j$ with $p_i=p_j$,
then, starting from all 0's, the probability of poisoning in states
$i$ or $j$ is 0. In particular, if all of the $p_i$'s are equal, 
then, starting from all 0's, the probability of poisoning is 0.
\end{theorem}

Theorem \ref{first theorem} is proved in Section 2,
Theorem \ref{second theorem} is proved in Section 3, 
some further comments and questions are given 
in Section 4 and finally there is an appendix which provides
some assistance for one of the proofs.

\section{\textbf{Poisoning in 1-dimension with rate $<1/2$}}
\noindent
In this section, we prove Theorem \ref{first theorem}.

\medskip\noindent
{\textit{Proof of Theorem 3:}}
The method used will be to consider a block of
adjacent 1's and find conditions under which they tend to spread.
At the right-hand end of such a block there must be a $0$. We
condition on the types occupying the next three spaces to the right.
There are 18 essentially different possibilities as the reader can check;
see column 1 of Table 1 for a list of these. Blocks such as
$030$ and $020$ or $203$ and $304$ are considered equivalent.

Consider the same model but defined only on the nonnegative integers.
Equivalently, at time 0, there is a 0 in position $-1$ and no particles
are allowed to arrive at this position. 
We show that if there is a $1$ at the origin,
then with a fixed positive probability, independent of the configuration
to the right of the origin, the origin remains in state $1$ forever and
the block of $1$'s containing the origin approaches $\infty$.
(Of course if the origin remains in state $1$ forever, it must be the case
that the block of $1$'s containing the origin approaches $\infty$.)
If this can be done, it follows that in the original model, any 1
will spread in both directions poisoning all of the integers with
a fixed positive probability. Since there are initially infinitely
many 0's or 1's, it easily follows that poisoning in state 1 occurs a.s.

In order to prove the statement concerning the spread
of 1's, for each configuration $\eta$
on the nonnegative integers which has a 1 or 0 at the origin, we
define its weight $W(\eta)$ as follows.
Let $B(\eta)$ be the block of 1's starting from the origin
in the configuration $\eta$. (This may be empty or infinite.) 
If $|B(\eta)|=\infty$ (equivalently $\eta\equiv 1$), then we let
$W(\eta):=\infty$. Otherwise, $W(\eta)$ is defined to be
$|B(\eta)|$ plus the score, as defined in column 2 of Table 1, 
of the block following the first 0 after $B(\eta)$. 
The idea is to define the scores of the blocks $(abc)$ in
such a way that the expected change in the weight is always positive, in
which case we will obtain a submartingale. 

\begin{center}
{\textbf{Table 1}}\\ {\textbf{Scores for essentially different
blocks}}

\begin{tabular}{|c|c|c|c|}\hline
block&Score&Expression&Follower\\\hline 
222&0.000&0.0007&2\\\hline
220&0.163&0.0012&2\\\hline 202&0.295&0.0022&2\\\hline
203&0.339&0.0002&3\\\hline 022&0.354&0.0034&2\\\hline
200&0.404&0.0002&2\\\hline 201&0.493&0.0018&00\\ \hline
020&0.498&0.0031&2\\\hline 002&0.570&0.0032&2\\\hline
000&0.664&0.0055&2\\\hline 001&0.827&0.0058&00\\\hline
010&0.920&0.0044&2\\\hline 102&1.008&0.0034&22\\\hline
100&1.157&0.0036&22\\\hline 011&1.173&0.0060&00\\\hline
101&1.456&0.0056&00,02\\\hline 110&1.555&0.0040&222\\\hline
111&1.997&0.0054&0222\\\hline

\end{tabular}
\end{center}

Letting $\{\eta_t\}_{t\ge 0}$ denote our process on the nonnegative integers,
a very long and tedious calculation, left to the reader, shows that 
there exists $c > 0$ such that for all 
$\eta$ with $B(\eta)\neq \emptyset$,
\begin{equation}\label{maineqn}
\frac{dE[W(\eta_t)|\eta_0=\eta]}{dt}|_{t=0} \ge c.
\end{equation}
[While this long detailed calculation 
is being left to the reader, the appendix contains a discussion which
aids the reader in making this calculation; perhaps it takes 2 hours of
work to check the above with the aid of the appendix and 5 hours otherwise.]
It follows that with positive probability, uniform in $\eta$,
with $B(\eta)\neq \emptyset$,
the block of 1's will grow to $\infty$ before the 1 at the origin
changes which is what we wanted to show. To carefully do this,
we follow the argument for Theorem 3 in Sudbury (1999) and proceed
as follows. Let $T$ be the (possibly infinite) stopping time
when the block of 1's at the origin disappears. Equation (\ref{maineqn})
implies that $\{W(\eta_{t\wedge T})\}_{t\ge 0}$ 
is a submartingale with respect to the natural filtration of
$\{\eta_t\}_{t\ge 0}$. Moreover, it can be shown that 
Equation (\ref{maineqn}) implies that
for $\epsilon$ sufficiently small, for any initial configuration $\eta$
with $B(\eta)\neq \emptyset$,
$U_t:=1-(1-\epsilon)^{W(\eta_{t\wedge T})}$
is a bounded submartingale and thus tends a.s.\ to a limit $U_\infty$. 
Let us assume further for the moment that $|B(\eta_0)| \ge 2$.
This gives us that $U_0\ge 1-(1-\epsilon)^2$. 
If the block of 1's dies out, then 
$U_\infty\le  1-(1-\epsilon)^{1.997}$. Since $E[U_0]=E[U_\infty]$,
it follows that it cannot be the case that 
the block of 1's dies out a.s. Since however $U_t$ converges
a.s., it must be the case
that the block of 1's grows to infinity with a uniform positive probability,
independent of $\eta_0$ with $|B(\eta_0)| \ge 2$.
If $|B(\eta_0)|=1$, it is clear that with a uniform positive probability,
independent of $\eta_0$, the single 1 at the origin spreads
and reaches size at least 2 at time 1. At that point, one can apply the
previous argument.
$\blacksquare$

\medskip\noindent
{\bf Remarks:} \\
(i). The score $(abc)$ represents whether this triple is
likely to aid the 1's in spreading or make them more likely to
contract. \\
(ii). The above proof is {\sl not} computer assisted in the sense that
one can check the correctness of the proof without a computer;
it suffices to use a hand calculator or in fact even hand calculations
suffice
(the latter requiring a good deal of patience). Nonetheless, a computer was
essential in helping us find the proof and in particular helping us
find the scores for the various blocks of length $3$. More
discussion concerning this point follows below.\\
(iii). One might hope that we could have carried out the above proof
using blocks of length 2 instead of blocks of length 3 but it seems that
this is not possible if we want to have $p_1 < 1/2$ as also
described below.

\medskip\noindent
While it is not needed at all for checking the correctness of the proof of
Theorem \ref{first theorem}, we now nonetheless
explain in some detail how the scores in Table 1 were arrived
at which in turn 
allowed us to obtain the proof. We first tried blocks of size 2.
Our method was to try various values of $p_1$. As in the proof, 
we considered only the situation at the right-hand end. We wrote down the
equations for the instantaneous rates of change of the expected
weight for the various blocks of size 2.
The basic idea was to find the values of the $(ab),ab
\neq 22$ which make all these rates 0. Having done this we 
wrote down the equation for $(22)$. If the rate was $>0$ then we
would have found a value of $p_1$ which allows the block of
1's to spread. If
the rate was $<0$, then we would need a larger value of $p_1$. (In
fact we found values of the $(ab)$ which gave a rate very close to
$0$ and then increased $p_1$ so that the rate of change of the
expected weight was strictly positive for all blocks).
When we computed this change in the expected weight,
we computed it assuming to the right of the block one had what one would
guess is the most disadvantageous scenario for
spreading of the 1's. The only case considered was $p_2=p_3=p_4$. A little
experimentation with other possibilities suggested they would give
less favorable results.

Unfortunately it was impossible to make all the necessary expressions
positive with $p_1<0.5$. This meant we had to go up to the next
level, looking at blocks of length 3  to the right of the block of
1's and their immediate 0. There were then 18 fundamentally
different expressions that must be made $>0$.  Again, we assumed that
$p_2=p_3=p_4$. We assumed 
that $1>0>2,3,4$ in the sense that the score for every block is
increased if the value in the block is increased by replacing it
by a ``higher'' value. For example, we assumed that
$(020)<(000)<(010)$. This assumption is easily checked by looking
at the 18 scores in column 2 of Table 1. Assuming this worse case scenario to
the right of the 3 block, one obtains the following rates of change.
The assumed worse case scenario to the right of a 3 block is listed in the
fourth column in Table 1.

\noindent
{\underline{000}}$$
p_1((002)+1+(100)+(010))
+3p_2(-1+(200)+(020)+(002)/3)+2p_2(000)-3(000)
$${\underline{203}}$$
p_1((002)+.5((002)+(200)))+ p_2(-1+(020)+(200)+(002)) $$ $$
+p_2(-1+(020)+(002)+((002)+(200))/2)-2(203)
$${\underline{100}}$$
p_1(2+(022)+(110)+(100))+3p_2((010)-1+(000))/2 $$ $$
+3p_2(000)+2p_2(100)+p_2(102)-3(100)
$${\underline{111}}$$
p_14+3p_2(-1+(011)2)/2-(111)+(1-p_1)((110)-(111))
$${\underline{001}}$$
p_1(1+(010)+(101)+(011))
+3p_2(-1+(000)+(201)+(000))-3(001)+(1-p_1)((000)-(001))
$${\underline{102}}$$
p_1(2+(100))+3p_2(-1+(010)+(002))/2 +p_2((002)+(100)+(002))-2(102)
$${\underline{220}}$$
p_1((020)+((200)+(220))/2)+p_2(-1+(022)) $$ $$
+p_2(-1+(022)+(020)+(220)+(200))-2(220)
$${\underline{022}}$$
p_1(1+(002))+3p_2(-1+(002))+2p_2(002)-2(022)
$${\underline{202}}$$
p_1((002)+((200)+(002))/2)+p_2(-1+(020))
+p_2(-1+(020)+(002)+(002)+(200))-2(202)
$${\underline{200}}$$
p_1((000)+(000)+(200))+p_2(-1+(020)+(220)+(202)) $$ $$
+p_2(-1+(020)+(000)3+2(200))-3(200)
$${\underline{020}}$$
 p_1(1+(202)+1.5(000)+.5(020))+p_2(-1+(002)+(220)+(022))
$$
$$
+2p_2(-1+(002)+1.5(000)+.5(020))-3(020)
$${\underline{201}}$$
p_1((001)2)+p_2(-1+(020)+(200)-1+(020)+2(001)+(200))-2(201)+(1-p_1)((200)-(201))
$${\underline{010}}$$
p_1((1+(102)+(110)+(010))
+3p_2(-1+(001)+(000))+p_2((000)2+(010))-3(010)
$${\underline{002}}$$
p_1(1+(022)+(102)+(000))+p_2(-1+(000)+(202)+(022))+2p_2(-1+(000)+(302)+(000)
)-3(002)
$${\underline{011}}$$
p_1(1+(110)+(111))+3p_2(-1+(001)+(001))-2(011)+(1-p_1)((010)-(011))
$${\underline{110}}$$
p_1(3+(110))+3p_2(-1+(011)+(010))/2+p_2(2(100)+(110))-2(110)
$${\underline{101}}$$
p_1(2+(102)+(111))+3p_2(-1+(010)+(001)+(001)+(100))/2-2(101)+(1-p_1)((100)-(101))
$${\underline{222}}$$
p_1(022)+p_2(-2+3(022))$$

A computer program was
written which calculated these expressions for blocks of any
length. The program goes through all possible sets of right-hand
blocks made up of $0,1,2,3,4$, finds the $0$'s and in turn
replaces them with $1,2,3,4$. The score of each block was then
replaced by the score that makes the relevant expression above
equal $0$. For example, the third to last equation above becomes

$$(110)_{new}=[p_1(3+(110))+3p_2(-1+(011)+(010))/2+p_2(2(100)+(110))]/2.$$

The first 17 equations were then iterated until the differences
between the new and old scores given to the blocks were
sufficiently close to 0. All 18 equations were not used because,
since the total rates in and out of the blocks must balance, the
set of equations would be singular. The last expression, which
is for $(222)$ was then tested. If it was positive, a suitable set
of values would have been found. If it was not,
 $p_1$ was replaced by a larger value and the procedure was started again.
A suitable set of scores was found for $p_1=0.4699$. The
expressions above were then independently checked (with the
package Minitab) using $p_1=0.47$ to ensure all the expressions
came out strictly positive. The first few decimal places of these expressions
are given in column 3 of Table 1. (It should be noted that the computer
program tests the transitions for all possible scenarios
to the right of the blocks of length 3 given above since
these scenarios can influence the outcome. The expressions given
above use simply the best guess as to which of the values
$0,1,2,3,4$ will give the least advantageous result for the spread
of the block of 1's. The guesses are in fact correct as can be
checked.) The checks made by minitab can then also be carried out by hand 
and/or pocket calculator obtaining the rigorous noncomputer assisted
proof presented above.

One can next attempt to prove that even lower values of $p_1$ can spread
with 4 types or alternatively study this problem with a different number of 
types. In either case, one has to deal with blocks of length greater than 
3 and then the equations become too complicated to exhibit or
to do calculator or hand computations with.
The computer program used can only
continue up until blocks of length 5 as the particular version of
Fortran used only allows arrays of size up to 4096.
Further, one can consider the case when there are infinitely many
gases. In this limiting case it is assumed that any new arrival
will react with an adjacent gas unless it is of type 1 in which
case there will be no reaction if the arriving gas is also type 1.
Since this situation is more favorable to gas 1, it is possible
to use only blocks of length 2 to show that in this case,
gas 1 can spread with a rate of 0.46.

Theorem \ref{first theorem} may then be improved as follows using the
computer program mentioned above.  Belief in the following theorem
thus relies on belief that the program is correct. The program 
gives the same results as hand calculation for blocks of length 2 and 3.
\begin{theorem} (Computer assisted proof) In one dimension, poisoning will 
occur for rates of gas 1 greater than those appearing in Table 2 
(with all other gases having the same rate) if the
initial state of the surface has an infinite number of 0's or 1's.
\end{theorem}
\begin{center}
{\textbf{Table 2}}\\ {\textbf{Upper bounds to critical values for
poisoning}} \hspace{1cm}
\begin{tabular}{|c|c|c|}\hline
Number of Gases&Rate of gas 1&Rate of other gases\\\hline
2&$>0.5$&$<0.5$\\\hline 3&0.445(0.40)&0.277\\\hline
4&0.432(0.38)&0.189\\\hline $\infty$&0.46(0.37)&0.000\\\hline

\end{tabular}
\end{center}
Figures in parentheses are the critical values suggested by
simulation. The case of 4 gases is a worse approximation to the
correct value than is 3 gases, as the computer program only allows the
use of blocks of 5 places in this case.

\section{\textbf{Nonpoisoning}}

In this section, we prove Theorem \ref{second theorem}.
The argument here is very similar to that used in
Bramson and Griffeath (1989) and Cox and Klenke (2000).

Recall that in this result, we are starting with all $0$'s.
For each $i=1,\ldots,n$, let $A_i$ be the event that poisoning occurs
and that the final state of all sites is $i$. Clearly, the $A_i$'s are disjoint
and their union is the event in question.

\begin{lemma} \label{01lemma}
For each $i$, $P(A_i) \in \{0,1\}$. (This lemma does not require any assumption on the
rates.)
\end{lemma}

\noindent
Given this lemma, we are done as follows. $P(A_i)=P(A_j)$ by symmetry and by
by Lemma \ref{01lemma}, these are each 0 or 1. Since they are disjoint, they each
must be 0. The last statement of the result follows immediately.
$\blacksquare$

\medskip\noindent
{\sl  Proof of Lemma \ref{01lemma}.}
Of course, we just need to consider $A_1$.
The idea is that although we have stochastic dynamics, we can encode all the randomness we will
need to drive the dynamics into random variables associated to each lattice point.
To drive the dynamics, for each location, we need $n$ independent
Poisson processes (for the arrivals of the $n$ different types of particles)
and we will also need some more random variables at each
location to be used to decide which particle will be reacted with
if the arriving particle lands next to more than 2 particles of a different type.

The details of how to do this are as follows.
Let $\{U_i\}_{i\in Z^d}$ be i.i.d.\ random variables each uniformly distributed on $[0,1]$.
The entire evolution of our process will simply be a function of the $\{U_i\}$'s
(i.e., it will be determined by the $\{U_i\}$'s); no further
randomness will be needed. On the unit interval $[0,1]$ with the Borel sets
and with Lebesgue measure,
define random variables $\{X_t\}_{t\ge 0}$, $\{Y_k\}_{k\ge 1}$ and $\{Z_k\}_{k\ge 1}$ such that
these three processes are independent of each other, the first process is a rate
1 Poisson process, the $Y_k$'s are i.i.d. with $P(Y_k=j)=p_j$ for $j=1,\ldots,n$,
and the $Z_k$'s are i.i.d. with each $Z_k$ being a uniform random ordering
of the $2d$ neighbors of the origin in $Z^d$. Note, crucially, that the
$X_t$'s, $Y_k$'s and $Z_k$'s are {\sl functions} defined on $[0,1]$.

Note that we know that on {\sl some} probability space, we can define random variables
with the above distribution but the point is that we want to define them on
the unit interval $[0,1]$ with the Borel sets. It is well
known that this probability space is
{\sl rich} enough to be able to define these random variables on it. We of course,
as usually done in probability theory, do not need to explicitly give what these functions
are; it is only required that they have the right distribution.

Now, for each $i\in Z^d$, we can consider the random variables
$\{X_t(U_i)\}_{t\ge 0}$, \\
$\{Y_k(U_i)\}_{k\ge 1}$ and $\{Z_k(U_i)\}_{k\ge 1}$. These
are independent for different $i$ (since the $U_i$'s are)
and have the same distribution as the
$\{X_t\}_{t\ge 0}$, $\{Y_k\}_{k\ge 1}$ and $\{Z_k\}_{k\ge 1}$ defined above
(because $U_i$ is uniform on $[0,1]$). Note that $X_t(U_i)$ (as well as these others)
are random variables since they are a composition of random variables. Note that
for this to hold, it was important that the underlying probability space for
$\{X_t\}_{t\ge 0}$, $\{Y_k\}_{k\ge 1}$ and $\{Z_k\}_{k\ge 1}$ was the unit interval
with the Borel sets and not the unit interval with the collection of Lebesgue measurable sets.

We now use these random variables to drive our dynamics as follows. Fix $i\in Z^d$.
The arrival times of particles at site $i$ will be taken to be
the Poisson process $\{X_t(U_i)\}_{t\ge 0}$.
The type of the $k$th arriving particle at site $i$ will be taken to be $Y_k(U_i)$.
If the $k$th particle arrives and there is at most 1 neighboring particle of a different type,
we of course know what to do. However, if there is more than 1 such neighboring particle,
then we look at the random ordering $Z_k(U_i)$ of the neighbors of the origin
and have our arriving particle choose, among those neighbors which have
a particle of a different type, the neighbor $i+u$ with the highest $u$ value according
to the ordering $Z_k(U_i)$, and then $i$ and this chosen site react and become 0. It is clear
that this generates our interacting particle system.

Now the event $A_1$ is of course measurable with respect to
$\{U_i\}_{i\in Z^d}$ and is also trivially translation invariant.
Since any i.i.d.\ process is ergodic (see for example, Walters (1975)),
it follows that $P(A_1) \in \{0,1\}$, as desired.
$\blacksquare$

\section{\textbf{Further remarks and conjectures}}

We believe the following strengthening of Theorem \ref{second theorem} should be true.

\begin{conjecture}  \label{C1}
Let $n$ and $d$ be arbitrary. If there are $i\neq j$ such that
$p_i=p_j$ and all the other $p_k$'s are no larger, then, starting from all 0's,
the probability of poisoning is 0.
\end{conjecture}

By Theorem \ref{second theorem},
the probability that we get poisoned in states $i$ or $j$ is 0. It seems natural
that it should be even harder to get poisoned in one of the other states with a smaller
$p_k$. The following monotonicity result seems reasonable and would imply
Conjecture \ref{C1}. Of course, this monotonicity result is not necessary for
Conjecture \ref{C1} and it cannot be ruled out that there are no other types of
monotonicity results which could be used instead to obtain Conjecture \ref{C1}.

\begin{conjecture}  \label{C2}
Consider two systems with the same $d$ and $n$ with vectors
$(p_1,p_2,\ldots,p_n)$ and
$(\tilde{p}_1,\tilde{p}_2,\tilde{p}_n)$ such that
$p_1 \le \tilde{p}_1$ and $p_i \ge \tilde{p}_i$ for $i\ge 2$.
Then the process of $1$'s in the first process is stochastically smaller than
the process of $1$'s in the second process.
\end{conjecture}

We point out however that the simplest naive way of trying to prove Conjecture \ref{C2}
doesn't work. Consider the case where $d=1$, $n=3$ and our two systems have rates
$(1/3,1/3,1/3)$ and $(1/2,1/4,1/4)$. The obvious way to couple these systems is to let
particles fall at the same location at the same time with the probability of the pair
$(i,j)$ falling being given by
$$
p_{(2,2)}=p_{(3,3)}=1/4,p_{(1,1)}=1/3,p_{(2,1)}=p_{(3,1)}=1/12
$$
and if both particles have a choice with whom to react with, they choose the same one.
Now, under this coupling, it can in fact happen that
there is a some location having a 1 in the first system but not having a 1 in the
second system. In the following realization, we first have
the first process getting a 2 and the second process getting a 1. After this,
in the next 5 steps, we have the same type particle arrive in the two 
systems and the types, in order of arrival, are 2,3,3,3, and 1. What we then see is

\medskip

\noindent
0000000000,0000000000

\noindent
0000200000,0000100000

\noindent
0000220000,0000000000

\noindent
0000220000,0000300000

\noindent
0000020000,0003300000

\noindent
0000000000,0003303000

\noindent
and finally

\noindent
0000010000,0003003000

\noindent
or

\noindent
0000010000,0003300000.

\section{\textbf{Appendix}}

In this section, we aid the reader in verifying Equation (\ref{maineqn})
by suggesting how this should be done. \\
(1). The reader should first check that the 18 expressions in Section 2 correspond
to the rate of change of the expected weight in what the reader would think
is the most disadvantageous scenario to the right of the block.
In column 4 of Table 1, these supposed most disadvantageous scenarios to the 
right of the block are listed. The case 101 is special in that
for certain terms, it is assumed that there is a 02 to the right and for
other terms, it is assumed that there is a 00 to the right. \\
(2). Next, one should plug in the scores for the blocks given in Table 1 into all the 
expressions for the blocks and check that all the numbers one obtains are 
positive. The first few decimals for these numbers are listed in column
3 of Table 1. \\
(3). Last, one has to consider all the other possible scenarios which can 
in fact sit to the right of the block, compute how each one affects 
the rate of change of the expected weight and check that all these 
other cases yield a larger value.  

One should observe that our assumed
worse case scenario sometimes holds uniformly and sometimes not in the
following sense. A 00 to the right of the block 001 is always worse
no matter what arrives next.
Similarly, a 2 to the right of the block 022 is always worse than a 0.
However, for the block 000, it is not uniformly worse to have a 2 to
the right rather than a 0. Should a 3 arrive
at the right most 0, it would have been worse to have a 0 to the right
since then the 3 would have remained but with a 2 to the right, the 3 would not
stay. However, if a 1 would have arrived at that position, it would be worse to
have a 2 since then the 1 would not stay. However, in all cases, when one
averages over all the possibilities, these assumed worse cases are in fact
worse case. However, this needs to be checked and this is precisely step (3) above.

\end{document}